# Mapping tori of free group automorphisms are coherent

By Mark Feighn and Michael Handel*

**Abstract**

The *mapping torus* of an endomorphism $\Phi$ of a group $G$ is the HNN-extension $G*_G$ with bonding maps the identity and $\Phi$. We show that a mapping torus of an injective free group endomorphism has the property that its finitely generated subgroups are finitely presented and, moreover, these subgroups are of finite type.

## 1. Introduction

A group is *coherent* if its finitely generated subgroups are finitely presented. Free groups are obviously coherent. The classification of surfaces and the fact that every cover of a surface is itself a surface, imply that surface groups are coherent. In the early 1970's, Scott [Sco73] and Shalen (unpublished) independently answered a question of Jaco by showing that the fundamental group of a 3-manifold is coherent. Stallings [Sta77] showed that $F_2 \times F_2$ is not coherent.

Since most finitely generated groups are not finitely presented, the question of which groups are coherent has centered on groups with special properties. Rips [Rip82] gave examples of incoherent small cancellation groups. Wise [Wis98] gave examples of compact negatively curved two complexes with incoherent fundamental group. McCammond and Wise [MW] have recently developed methods that allow them to show, for example, that a one-relator group $\langle A|W^n \rangle$ is coherent for all large $n$.

Many authors have explored parallels between the mapping class group $\mathrm{MCG}(S)$ of a compact surface and the outer automorphism group $\mathrm{Out}(F_n)$ of the free group on $n$ letters. Various analogues of Thurston's classification have been produced for $\mathrm{Out}(F_n)$ ([BH92], [CV86], [Lus92], [Sel96], [Sel]). The Scott conjecture (proven in [BH92]) bounds the rank of the fixed subgroup

*Both authors gratefully acknowledge the support of the National Science Foundation.

of a free group automorphism and is a generalization of a result of Nielsen for surface automorphisms. The $\text{Out}(F_n)$-analogue of the Nielsen problem for $\text{MCG}(S)$ (lift subgroups to the homeomorphism group of $S$) is to lift subgroups to the group of homotopy equivalences, up to a certain natural equivalence, of a marked graph. Results along these lines have proved by Culler [Cul84] and in [BFH96]. The Tits Alternative holds for $\text{MCG}(S)$ ([Iva84], [McC85]) as well as for $\text{Out}(F_n)$ ([BFH98], [BFH96]).

For a topological space $X$ and a map $f : X \to X$, the *mapping torus of $f$* is the quotient

$$M(f) := (X \times I)/\sim$$

where $\sim$ is the equivalence relation generated by $(f(x), 0) \sim (x, 1)$. Mapping tori of surface automorphisms have played a significant role in the study of 3-manifolds. For example, a part of Thurston's hyperbolization theorem is that the mapping torus of a pseudo-Anosov automorphism of a compact surface is hyperbolic [Thu]. In fact, it is conjectured by Thurston that every finite volume hyperbolic 3-manifold is finitely covered by such a mapping torus.

For a group endomorphism $\Phi : G \to G$, the *mapping torus*, denoted $M(\Phi)$, is the HNN-extension $G*_G$ where the bonding maps are the identity and $\Phi$. If the homomorphism $\Phi : \pi_1(X) \to \pi_1(X)$ is induced by $f : (X, *) \to (X, *)$ then $M(\Phi)$ is isomorphic to the fundamental group of $M(f)$. If $\Phi$ is an injective endomorphism then $M(\Phi)$ is also called an *ascending* HNN-*extension*.

It is natural to ask whether $M(\Phi)$ shares interesting properties with the class of 3-manifold groups. For example, if $\Phi$ is a hyperbolic automorphism of a word hyperbolic group (such as a finitely generated free group) then $M(\Phi)$ is word hyperbolic [BF92]. The main result of this paper is along this line.

THEOREM 1.1 (Main Theorem). *The mapping torus $M(\Phi)$ of an injective endomorphism $\Phi$ of a (possibly infinite rank) free group is coherent.*

This answers Problem 17 of Baumslag's 1973 problem list [Bau74]. Wise [Wis98] produces incoherent groups with presentations strikingly similar to that of $M(\Psi)$. Gersten [Ger81] showed the incoherence of the double $F_n *_H F_n$ of a free group $F_n$ of rank $n > 1$ over a finite index subgroup $H$ such that $[F_n : H] > 2$. Since $F_n *_H F_n$ is a subgroup of $F_n *_H$, this group is also not coherent. So, the theorem is, in a certain regard, sharp.

To give more detailed information about presentations of finitely generated subgroups of mapping tori, it is convenient to slightly generalize the notion of a mapping torus.

For a pair of topological spaces $Y \subset X$ and a map $f : Y \to X$, the *mapping torus of $f$* is the quotient

$$M(f) := (X \cup_{Y = Y \times \{0\}} (Y \times I))/\sim,$$

where $\sim$ is the equivalence relation generated by $f(y) \sim (y,1)$ for $y \in Y$. For a pair of groups $F \subset G$ and a homomorphism $\Psi : F \to G$, the *mapping torus*, denoted $M(\Psi)$, is the HNN-extension $G*_F$ where the bonding maps are the inclusion of $F$ into $G$ and $\Psi$. If the homomorphism $\Psi : \pi_1(Y) \to \pi_1(X)$ is induced by $f : (Y,*) \to (X,*)$ then $M(\Psi)$ is isomorphic to the fundamental group of $M(f)$. In this paper, we restrict attention to the class $\mathcal{M}$ of mapping tori where $G$ is a (possibly infinite rank) free group, $F$ is a free factor, and $\Psi$ is injective. Let $\mathcal{M}_0$ denote the subclass where additionally $G$ is finitely generated. By taking $F$ to be trivial, we see that $\mathcal{M}_0$ contains the class of nontrivial, finitely generated free groups. It is easy to check (see Proposition 2.1) that every element of $\mathcal{M}$ can be realized as the mapping torus of an injective endomorphism of a (possibly infinite rank) free group. Thus the class $\mathcal{M}$ is the same as the class of mapping tori of injective free group endomorphisms. The class $\mathcal{M}_0$ properly contains the class of mapping tori of injective endomorphisms of finitely generated free groups.

Each element $M(\Psi) \in \mathcal{M}$ has a class of *preferred presentations*

$$\langle t, A, B, | C \rangle$$

where $A$ is a basis for $F$, $A \amalg B$ is a basis for $G$, and $C = \{tat^{-1}(\Psi(a))^{-1} | a \in A\}$. For an element of $\mathcal{M}_0$, we also require that in a preferred presentation the cardinalities of the sets $A$, $B$, (and hence $C$) be finite. A group is of *finite type* if it has a compact Eilenberg-Mac Lane space. The 2-complex associated to a preferred presentation of an element of $\mathcal{M}$ is an Eilenberg-Mac Lane space (see [SW79]). So, groups in $\mathcal{M}_0$ are of finite type.

We can now give a more precise statement of our main result.

THEOREM 1.2.   *A nontrivial finitely generated subgroup of an element $M(\Psi) \in \mathcal{M}$ is in $\mathcal{M}_0$. In particular, each finitely generated subgroup of $M(\Psi)$ has finite type and $M(\Psi)$ is coherent.*

Our methods are geometric. In fact, our main technique (see §4) may be viewed as a relative version of Stallings folds [Sta83] (see §3).

We thank Dani Wise for pointing us to the Baumslag problem list and for providing us relevant examples of incoherent groups. We learned of this problem from Peter Shalen whom we thank for a very interesting history of the 3-manifold result.

## 2. Preliminaries

Throughout this paper, we will use the following notation. The free group with (not necessarily finite) basis $\mathcal{E} = \{e_i | i \in I\}$ is denoted $\mathbb{F}$, and $\Phi$ will denote an injective endomorphism of $\mathbb{F}$. As described in the introduction, the *mapping torus* $M(\Phi)$ is the HNN-extension $\mathbb{F}*_{\mathbb{F}}$ with bonding maps the identity and $\Phi$. By definition, $M(\Phi)$ has the presentation with generators $\mathcal{E} \cup \{t\}$ and relations $\{te_it^{-1}(\Phi(e_i))^{-1}|i \in I\}$. (Here and throughout we slightly abuse notation and write $\Phi(e_i)$ for both the element of $\mathbb{F}$ and the freely reduced word in $\mathcal{E}$ representing it.) Since $M(\Phi)$ is an HNN-extension with injective bonding maps, $\mathcal{E}$ freely generates a subgroup of $M(\Phi)$ that we identify with $\mathbb{F}$ (see, for example, [Ser80, Cor. 1, p. 45]). We denote by $p: M(\Phi) \to \mathbb{Z}$ the homomorphism defined by $p(t) = 1$ and $p(e_i) = 0$ for $i \in I$.

It is well-known that an element of $\mathcal{M}$ is also a mapping torus of an injective free group endomorphism. For completeness, we include a proof.

PROPOSITION 2.1. *Let $F$ be a free factor of $\mathbb{F}$ and $\Psi: F \to \mathbb{F}$ be an injective homomorphism. There is a free group $\mathbb{F}'$ and an injective endomorphism $\Phi: \mathbb{F}' \to \mathbb{F}'$ such that $M(\Phi)$ is isomorphic to $M(\Psi)$.*

*Proof.* We may assume that $\{e_j|j \in J\}$ is a basis for $F$ where $J \subset I$. Take $\mathbb{F}'$ to have the basis $\{e_{i,0}|i \in I\} \cup \{e_{i,k}|i \in I \setminus J, k = 1, 2, \cdots\}$. Denote by $\mathbb{F}'_1, \mathbb{F}'_2$, and $\mathbb{F}'_3$ the free factors of $\mathbb{F}'$ generated by, respectively, $\{e_{j,0}: j \in J\}$, $\{e_{i,0}: i \in I \setminus J\}$, and $\{e_{i,k}: i \in I \setminus J, k \geq 1\}$; thus $\mathbb{F}' = \mathbb{F}'_1 * \mathbb{F}'_2 * \mathbb{F}'_3$. After identifying each $e_{i,0} \in \mathbb{F}'$ with $e_i \in \mathbb{F}$, we may view $\Psi$ as a map from $\mathbb{F}'_1$ to $\mathbb{F}'_1 * \mathbb{F}'_2$. Define $T: \mathbb{F}'_2 * \mathbb{F}'_3 \to \mathbb{F}'_3$ by $T(e_{i,k}) = e_{i,k+1}$ and define $\Phi = \Psi * T : \mathbb{F}' \to \mathbb{F}'$. Since $\Psi$ is an injective endomorphism and $T$ is an isomorphism, $\Phi$ is an injective endomorphism.

By construction, $e_i \mapsto e_{i,0}$ and $t \mapsto t$ defines a homomorphism from $M(\Psi)$ to $M(\Phi)$ and $e_{i,k} \mapsto t^k e_i t^{-k}$ and $t \mapsto t$ defines a homomorphism from $M(\Phi)$ to $M(\Psi)$. These homomorphisms are inverses and hence isomorphisms. $\square$

We gather some elementary observations in the following lemma.

LEMMA 2.2.
(1) *Every element $g \in M(\Phi)$ has a representation of the form $g = t^{-q}xt^r$ where $q, r \geq 0$ and $x \in \mathbb{F}$.*

(2) *If a finitely generated subgroup $H$ of $M(\Phi)$ contains $t$ then there is a finite set $A \subset \mathbb{F}$ such that $H = \langle t, A \rangle$.*

(3) *If a subgroup $H$ of $M(\Phi)$ contains $t$ then $\Phi(H \cap \mathbb{F}) \subset H \cap \mathbb{F}$.*

(4) *If $g \in M(\Phi)$ satisfies $p(g) = m > 0$, then $g \in \langle \mathbb{F}, t^m \rangle$.*

*Proof.* Since $tx = \Phi(x)t$ and $xt^{-1} = t^{-1}\Phi(x)$ for all $x \in \mathbb{F}$, given any word in $\{t\} \cup \mathcal{E}$, we can move positive powers of $t$ to the right and negative powers of $t$ to the left. This implies (1) which in turn implies (2). Item (3) follows from the fact that $\Phi(x) = txt^{-1}$ for all $x \in \mathbb{F}$. It remains to prove (4). If $p(g) = m > 0$, then (1) implies that $g = t^{-q}xt^{q+m}$ for some $x \in \mathbb{F}$ and some $q \geq 0$. Choose $j \geq 0$ so that $q + j$ is a multiple of $m$ and write $g = t^{-q}t^{-j}t^jxt^{q+m} = t^{-(q+j)}\Phi^j(x)t^{q+m+j} \in \langle \mathbb{F}, t^m \rangle$. $\square$

The following proposition is the heart of Theorem 1.2. Its proof occupies all of the remaining sections.

PROPOSITION 2.3 (Main Proposition).   *If $\Phi$ is an injective endomorphism of $\mathbb{F}$, and if $H$ is a finitely generated subgroup of $M(\Phi)$ that contains $t$, then $H \in \mathcal{M}_0$. In fact, $H$ has a preferred presentation of the form $\langle t, A, B | C \rangle$ where*

- $A = \{a_1, \cdots, a_m\}$, $B = \{b_1, \cdots, b_r\}$, and $C = \{r_1, \cdots, r_m\}$ are finite sets in $\mathbb{F}$,
- $r_j = t a_j t^{-1} w_j^{-1}$ for $w_j = \Phi(a_j)$ and $1 \leq j \leq m$, and
- $\langle A, \Phi(A) \rangle = \langle A, B \rangle$.

We conclude this section by reducing our main theorem to our main proposition.

*Proof of Theorem* 1.2 *assuming Proposition* 2.3. By Proposition 2.1, we may assume that $\Psi$ is an injective endomorphism of $\mathbb{F}$, which we call $\Phi$ to conform to the notation of Proposition 2.3.

Let $H$ be a finitely generated subgroup of $M(\Phi)$ and suppose, at first, that $H \subset \operatorname{Ker}(p)$. Lemma 2.2(1) implies that each $g \in H$ has the form $t^{-q} x t^q$ where $x \in \mathbb{F}$ and where $q \geq 0$. Let $i_{t^k}(g) = t^k g t^{-k}$ be the inner automorphism of $M(\Phi)$ determined by $t^k$. If $k \geq q$, then $i_{t^k}(g) = t^{k-q} x t^{q-k} = \Phi^{k-q}(x) \in \mathbb{F}$. Since $H$ is finitely generated, there exists $k > 0$ such that $i_{t^k}(H) \subset \mathbb{F}$ and so $H$ is free.

If $t \in H$, then $H \in \mathcal{M}_0$ by the Main Proposition and there is nothing more to prove.

For the general case we may assume that $p(H)$ is generated by some $m > 0$. Lemma 2.2(4) implies that $H$ is contained in the subgroup $M^* = \langle t^m, e_i \in \mathcal{E} \rangle$. Choose $g_m \in H$ such that $p(g_m) = m$. Lemma 2.2(1) implies that $g_m = t^{-p} b t^{m+p}$ for some $p \geq 0$ and $b \in \mathbb{F}$. Up to changing $H$ by an isomorphism, we may replace $H$ by $i_{t^p}(H)$ and $g_m$ by $bt^m$. Define $\Theta = i_b \Phi^m$. The assignments $t^m \mapsto b^{-1} s$ and $e_i \mapsto e_i$ for $e_i \in \mathcal{E}$ define an isomorphism between $M^*$ and $M(\Theta) = \langle s, e_i \in \mathcal{E} | s e_i s^{-1} = \Theta(e_i) \rangle$ that carries $H$ to a subgroup $H'$ containing $s$. The argument of the previous paragraph implies that $H'$, and hence $H$, is in $\mathcal{M}_0$. □

## 3. Labeled graphs and free groups

In this section, we recall a procedure of Stallings ([Sta83, Algorithm 5.4]) that from a finite set of words in the generators $\mathcal{E}$ of $\mathbb{F}$ produces a basis for the subgroup of $\mathbb{F}$ generated by these words. All of the results in this section are contained in [Sta83]. In Section 4, a relative version of this procedure will be presented.

*Definition* 3.1.   A graph is a one-dimensional CW-complex. The rose $R$ associated to $\mathcal{E}$ is the graph with one vertex $v_R$ and with oriented edges in one-to-one correspondence with $\mathcal{E}$. We identify $\pi_1(R, v_R)$ with $\mathbb{F}$ in the usual way: the homotopy class of any orientation-preserving immersion of $[0, 1]$ onto the edge of $R$ corresponding to $e_i$ is identified with $e_i$. A *labeled graph* $X$ is a finite connected graph with a basepoint $*$ and with oriented edges, each labeled by some $e_i \in \mathcal{E}$. The labeling defines, up to homotopy relative to the vertices of $X$, a map $f_X : X \to R$ that is injective when restricted to the interior of any edge. We write $\pi_1(X)$ for $\pi_1(X, *)$ and denote $f_\#(\pi_1(X)) \subset \mathbb{F}$ by $X^\#$. If $f_X$ is an immersion then we say that $X$ is *tight*. If $H$ is a subgroup of $\mathbb{F}$ and $X^\# = H$, then we say that $X$ is a *labeled graph for* $H$.

The importance of tightness is indicated by the following proposition ([Sta83, Prop. 5.3]).

PROPOSITION 3.2.   *If $X$ is a tight labeled graph then $(f_X)_\# : \pi_1(X) \to \mathbb{F}$ is injective.*

*Definition* 3.3 ([Sta83, §3]).   If $X$ is a labeled graph for $H$ which is not tight, then there is (at least one) pair of distinct closed edges $E_1$ and $E_2$ of $X$ with the same initial or terminal endpoints and the same labels. Let $q : X \to X'$ be the quotient map obtained by identifying $E_1$ with $E_2$, and let $q(*)$ be the basepoint for $X'$. The labeling for $X$ descends to a labeling of $X'$ such that $f_X = f_{X'} q$. We say that $X'$ is obtained from $X$ by *folding* $E_1$ and $E_2$, and we call the map $q$ a *fold*. A fold $q$ is a homotopy equivalence unless $E_1 \cup E_2$ is a *bigon*, i.e. they share both initial and terminal endpoints. In either case, $q_\#$ is onto and so $X'$ is also a labeled graph for $H$.

When a fold $q : X \to X'$ is not a homotopy equivalence, the rank of $X'$ is less than the rank of $X$ where by the *rank* of a graph, we mean its first betti number. In particular, folding does not increase rank. There are always fewer edges in $X'$ than in $X$. The next proposition follows easily.

PROPOSITION 3.4 ([Sta83, §3]).   *There is a finite sequence of folds that from a labeled graph $X$ for $H$ produces a tight labeled graph $\hat{X}$ for $H$ with $\mathrm{rank}(\hat{X}) \leq \mathrm{rank}(X)$. If $f_\#$ is not injective, then $\mathrm{rank}(\hat{X}) < \mathrm{rank}(X)$.*

*Notation* 3.5.   Suppose that $W = \{w_1, \ldots, w_n\}$ is a set of words in $\mathcal{E}$. For each $i$, let $X(\{w_i\})$ be the labeled graph for $\langle w_i \rangle$ that is homeomorphic to a circle and such that $f_{X(\{w_i\})}$ is an immersion away from the basepoint. Let $X(W)$ be the labeled graph $\bigvee_{i=1}^n X(\{w_i\})$ for $\langle W \rangle$ obtained by wedging at the basepoints.

We can now present Stallings' algorithm. Given a finite set $W$ of words in $\mathcal{E}$, let $\hat{X}(W)$ be the tight labeled graph obtained from $X(W)$ as in Proposition 3.4. By Proposition 3.2 $(f_{\hat{X}(W)})_\#$ is injective. Thus, if $A$ is a basis for $\pi_1(\hat{X}(W))$ then $(f_{\hat{X}(W)})_\#(A)$ is a basis for $\langle W \rangle$.

*Example* 3.6. Suppose that $\mathbb{F} = \langle e_1, e_2, e_3 \rangle$ has rank 3 and that

$$W = \{e_2 e_1 e_3, e_2 e_3 e_1, e_3^{-1} e_2 e_1, e_2 e_3 e_2^{-1} e_3\}.$$

The procedure outlined above will be used to find a basis for $\langle W \rangle$. Start with $X(W)$ which is depicted in Figure 1 below. In the figures, the direction of the arrows on an edge indicates the orientation of the edge, an edge with $i$ arrows is labeled $e_i$, and the large vertices are basepoints. To go from Figure 1 to Figure 2, two edges labeled $e_3$ are folded. From Figure 2 to Figure 3, four folds are performed. (In our diagrams, edges that are to be folded need not be particularly close. For example, consider the two edges labeled '$e_3$' incident to the base point in Figure 2.) Finally, from Figure 3 to the tight Figure 4, there are two folds. Choosing the edges in Figure 4 containing the basepoint as a maximal tree, the basis $\{e_2 e_1 e_3, e_3^{-1} e_2 e_1, e_2 e_3 e_1\}$ for $\langle W \rangle$ is obtained.

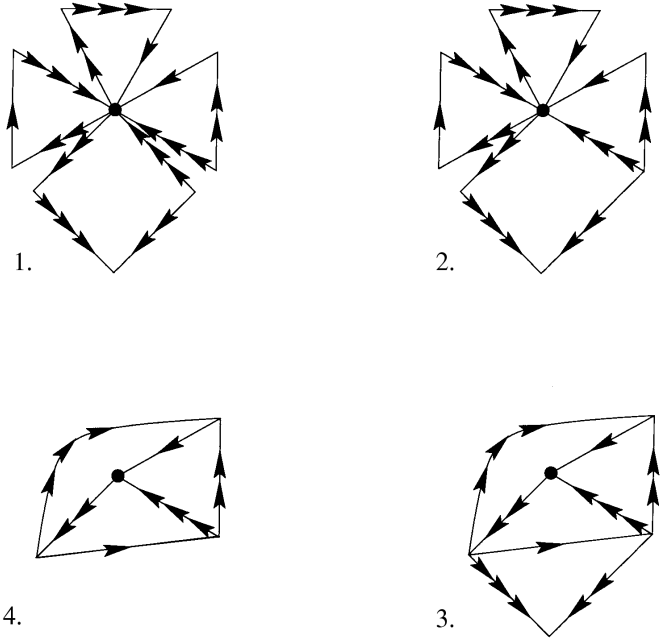

Example 3.6

## 4. Labeled graph pairs

Let $\Phi$ be an injective endomorphism of $\mathbb{F}$ and suppose that $H$ is a finitely generated subgroup of $M(\Phi)$ that contains $t$. To prove Proposition 2.3, we will need a mechanism for keeping track of pairs of subgroups of $\mathbb{F}$ of the form $(\langle A, B \rangle, \langle A \rangle)$. This leads us to consider pairs of labeled graphs.

*Definition* 4.1. A *labeled graph pair* is a pair $(Z, X)$ such that $X \subset Z$ are labeled graphs with a common basepoint. The pair $(Z, X)$ is *tight* if $Z$ is tight. The relative rank of $(Z, X)$ is defined by $\mathrm{rr}(Z, X) = \mathrm{rank}(\pi_1(Z)) - \mathrm{rank}(\pi_1(X)) = \chi(X) - \chi(Z)$.

If $H$ is a finitely generated subgroup of $M(\Phi)$ then a labeled graph pair $(Z, X)$ is a *labeled graph pair for $H$* if $\langle t, X^{\#} \rangle = H$ and $Z^{\#} \subset \langle X^{\#}, \Phi(X^{\#}) \rangle \subset H \cap \mathbb{F}$. (The last inclusion is a consequence of Lemma 2.2(3) and is included here for emphasis.) If additionally $Z^{\#} = \langle X^{\#}, \Phi(X^{\#}) \rangle$, or equivalently $\Phi(X^{\#}) \subset Z^{\#}$, then $(Z, X)$ satisfies the *invariance property* and $(Z, X)$ is an *invariant labeled graph pair for $H$*.

*Remark* 4.2. If $H$ is a finitely generated subgroup of $M(\Phi)$ that contains $t$, then, by Lemma 2.2(2), there is a finite set $A$ of $\mathbb{F}$ such that $H = \langle t, A \rangle$. If we take $X = X(A)$ and $Z = X(A \cup \Phi(A))$, then $(Z, X)$ is an invariant labeled graph pair for $H$.

The rest of this section is devoted to the details of a relative version of the first part of Proposition 3.4; the second part of Proposition 3.4 will be addressed in Section 5. The relative procedure will start with an invariant labeled graph pair for $H$ and produce a tight invariant labeled graph pair for $H$.

If $(Z, X)$ is a labeled graph pair then a fold $q_Z$ of $Z$ induces a map of pairs $q = (q_Z, \check{q}_X) : (Z, X) \to (Z_1, X_1)$ where $X_1 = q_Z(X)$. We still call $q$ a *fold*. Mimicking the absolute case, we will use folding to improve an invariant labeled graph pair $(Z, X)$ that isn't tight. We chose the notation $\check{q}_X : X \to X_1$ for the induced map to remind the reader of an important point—$\check{q}_X$ *need not be a fold of $X$*. In fact, it may be that $\mathrm{rank}(X_1) > \mathrm{rank}(X)$, perhaps resulting in a loss of the invariance property. In the next lemma, we record the possibilities for $\check{q}_X$. The proof, an immediate consequence of the definitions, is omitted.

LEMMA 4.3. *Suppose that $(Z, X)$ is a labeled graph pair, that $q : (Z, X) \to (Z_1, X_1)$ is a fold, and that $q_Z$ folds edges $E_1$ and $E_2$ of $Z$. Let $v$ denote a vertex shared by $E_1$ and $E_2$.*

(1) *If $E_1$ and $E_2$ are edges of $X$ then $\check{q}_X$ is a fold.*

(2) If $E_1 \cup E_2$ is not a bigon and is not contained in $X$, and if the vertices of $E_1$ and $E_2$ that are opposite $v$ are contained in $X$, then $\check{q}_X$ is defined by identifying a pair of distinct vertices in $X$.

(3) In all remaining cases $\check{q}_X$ is a homeomorphism.

*Notation* 4.4. In Case (1) we say that $q$ is a *subgraph fold*, and in Case (2) we say that $q$ is *exceptional*.

LEMMA 4.5. *Suppose that $(Z, X)$ is a labeled graph pair and that $q : (Z, X) \to (Z_1, X_1)$ is a fold. Then the following hold*:

(1) $Z^\# = Z_1^\#$.

(2) $X^\# \subset X_1^\#$.

(3) If $q$ is not exceptional and if $(Z, X)$ has the invariance property, then $(Z_1, X_1)$ has the invariance property.

(4) $\mathrm{rr}(Z_1, X_1) \leq \mathrm{rr}(Z, X)$.

(5) If $E_1 \cup E_2$ is a bigon that is not contained in $X$ then $\mathrm{rr}(Z_1, X_1) < \mathrm{rr}(Z, X)$.

(6) If $q$ is exceptional, then $\mathrm{rr}(Z_1, X_1) < \mathrm{rr}(Z, X)$.

*Proof.* Item (1) follows from $f_Z = f_{Z_1} q_Z$ and the fact that $(q_Z)_\#$ is onto. Item (2) follows from $f_X = f_{X_1}(\check{q}_X)$. If $q$ is not exceptional, then Lemma 4.3 implies that $(\check{q}_X)_\#$ is onto and hence that the inclusion in item (2) is an equality. Along with (1), this implies (3).

Suppose that $E_1 \cup E_2$ is a bigon. Then $\mathrm{rank}(Z_1) = \mathrm{rank}(Z)-1$, $\mathrm{rank}(X_1) = \mathrm{rank}(X) - 1$ if $E_1 \cup E_2 \subset X$, and $\mathrm{rank}(X_1) = \mathrm{rank}(X)$ otherwise. Suppose that $E_1 \cup E_2$ is not a bigon. Then $\mathrm{rank}(Z_1) = \mathrm{rank}(Z)$, $\mathrm{rank}(X_1) = \mathrm{rank}(X)$ if $q$ is not exceptional, and $\mathrm{rank}(X_1) = \mathrm{rank}(X) + 1$ if $q$ is exceptional. This implies (4), (5), and (6). □

Since an exceptional fold need not preserve the invariance property, we introduce a second move.

*Definition* 4.6. Suppose that $(Z, X)$ is an invariant labeled graph pair for $H$ and that $q : (Z, X) \to (Z_1, X_1)$ is a fold. We define a new invariant labeled graph pair $(Z_2, X_2)$ for $H$ as follows. There are three cases to consider; in all cases $X_2 = X_1$. If $q$ is not exceptional, then define $Z_2 = Z_1$. If $q$ is exceptional, let $p_1$ and $p_2$ be the points that are identified by $\check{q}_X$. Choose immersed paths $d_1$ and $d_2$ in $X$ from $*$ to $p_1$ and $p_2$ respectively, and let $\delta$ be the element of $X_1^\# \subset Z^\# \subset H$ represented by $f_{X_1}(\check{q}_X(d_1 d_2^{-1}))$. If $\Phi(\delta) \in Z_1^\#$, then define $Z_2 = Z_1$. Finally, if $\Phi(\delta) \notin Z_1^\#$, define $Z_2 = Z_1 \vee X(\{\Phi(\delta)\})$. We say that $(Z_2, X_2)$ is obtained from $(Z, X)$ by *folding and adding a loop if necessary*.

LEMMA 4.7. *Using the notation of Definition* 4.6, *the following hold.*

(1) $(Z_2, X_2)$ *is an invariant labeled graph pair for* $H$.

(2) *The number of vertices in $X_2$ is at most the number of vertices in $X$, and if the fold is exceptional, then $X_2$ has fewer vertices than $X$.*

(3) *If the fold is not exceptional, then $Z_2$ has fewer edges than $Z$.*

(4) $\mathrm{rr}(Z_2, X_2) \leq \mathrm{rr}(Z, X)$.

(5) *If $E_1 \cup E_2$ is a bigon that is not contained in $X$, then $\mathrm{rr}(Z_2, X_2) < \mathrm{rr}(Z, X)$.*

(6) *If the fold is exceptional and if $\Phi(\delta) \in Z_1^\#$ then $\mathrm{rr}(Z_2, X_2) < rr(Z, X)$.*

*Proof.* If the fold is not exceptional then $(Z_2, X_2) = (Z_1, X_1)$ and (1) follows from Lemma 4.5(3). If the fold is exceptional, then $Z_2^\# = \langle Z_1^\#, \Phi(\delta) \rangle = \langle Z^\#, \Phi(\delta) \rangle$ and $X_2^\# = X_1^\# = \langle X^\#, \delta \rangle$. In this case, (1) follows from Lemma 2.2(3) and the assumption that $(Z, X)$ is an invariant labeled graph pair for $H$.

Item (2) follows from Lemma 4.3 and the fact that $X_2 = X_1$. For a non-exceptional fold, $Z_2 = Z_1$. Item (3) follows immediately and the non-exceptional case of (4) follows from Lemma 4.5(4). The exceptional case of (4) follows from Lemma 4.5(6) and $\mathrm{rr}(Z_2, X_2) = \mathrm{rr}(Z_1, X_1) + 1$. If $E_1 \cup E_2$ is a bigon and $E_1$, say, is not contained in $X$, then $X_2 = X_1 = X$ and $Z_2 = Z_1$ is obtained from $Z$ by erasing the interior of $E_1$. This implies (5). Item (6) follows from Lemma 4.5(6) and the fact that $(Z_2, X_2) = (Z_1, X_1)$ in this case. $\square$

We can now define our relative version of the tightening procedure of Stallings.

*Definition* 4.8. Assume that $(Z, X)$ is an invariant labeled graph pair for $H$. If $(Z, X)$ is tight then do nothing. If $f_X$ is not an immersion, then perform a subgraph fold (Notation 4.4). If $f_X$ is an immersion and if $f_Z$ is not an immersion, then fold and add a loop if necessary. Repeat until the resulting invariant labeled graph pair for $H$, denoted $(\hat{Z}, \hat{X})$, is tight. Items (2) and (3) of Lemma 4.7 guarantee that this procedure terminates in finite time. We say that $(\hat{Z}, \hat{X})$ is obtained from $(Z, X)$ by *tightening*.

*Example* 4.9. Suppose that $\mathbb{F} = \langle e_1, e_2, e_3 \rangle$ has rank 3. Let $\Phi$ be the automorphism of $\mathbb{F}$ given by $\Phi(e_1) = e_2$, $\Phi(e_2) = (e_2^{-1}e_3e_2)$, and $\Phi(e_3) = e_2e_1^{-1}e_2$. Set $A = \{e_3^{-1}e_1, e_2^{-1}e_3^{-1}e_1^2e_3^{-1}e_1\}$. We have $\Phi(A) = \{e_2^{-1}e_1, e_2^{-1}e_3^{-1}e_1^2\}$. Set $Z = X(A \cup \Phi(A))$ and $X = X(A)$. Then $(Z, X)$ is an invariant labeled graph pair for the subgroup $= \langle t, A \rangle$ of $M(\Phi)$, and is pictured in Figure 1 below. We follow the conventions of the previous example, but add that edges of $Z$ that are not in $X$ are labeled with arrows that are not filled-in. The tightening procedure applied to $(Z, X)$ is pictured below. We refer to the closure of the complement of the subgraph as the overgraph. Figure 2 is obtained from Figure 1 by three subgraph folds. The subgraph in Figure 2 is tight. To go to Figure 3, a loop of the overgraph is folded into the subgraph. To go to Figure 4, part of the remaining loop of the overgraph is folded into the subgraph. To go to Figure 5, the edge of the overgraph is folded into the subgraph. This fold is exceptional and so we must add a loop labeled $\Phi(e_2^{-1}e_1) = e_2^{-1}e_3^{-1}e_2e_2$ which is depicted in Figure 6. Finally, Figure 7 shows $(\hat{Z}, \hat{X})$ and is obtained from Figure 6 by folding three edges of the overgraph into the subgraph. Note that $\mathrm{rr}(Z, X) = 2$ whereas $\mathrm{rr}(\hat{Z}, \hat{X}) = 1$.

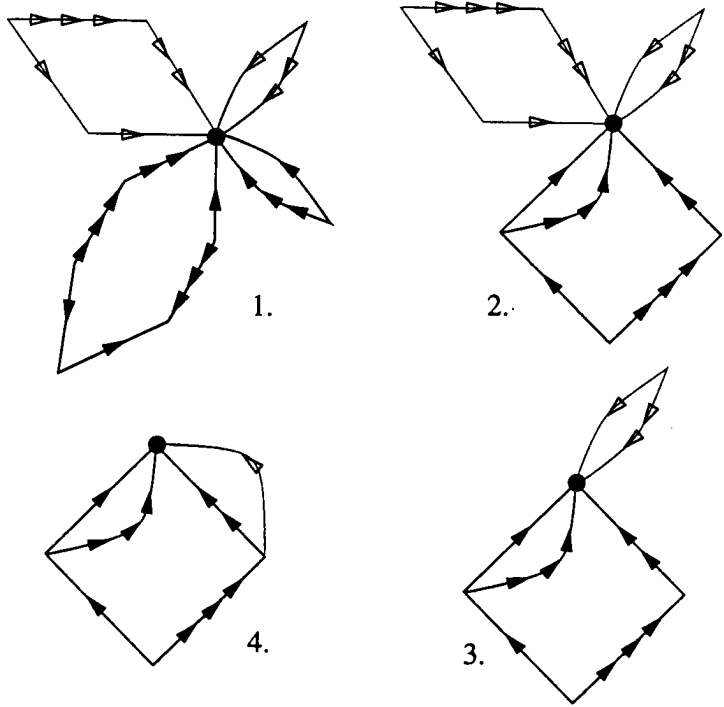

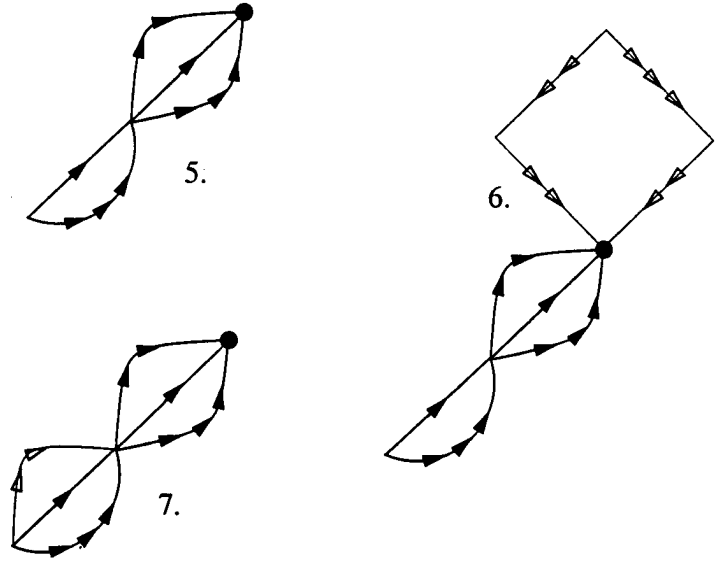

Example 4.9

## 5. Reducing relative rank

The goal of this section is a relative version of Proposition 3.4.

*Notation* 5.1.

- If $(Z, X)$ is a labeled graph pair, then the normal closure of the included image of $\operatorname{Ker}((f_X)_\#)$ in $\pi_1(Z)$ is denoted by $N(Z, X)$.

- A composition $Q = (Q_Z, Q_X) : (Z, X) \to (Z', X')$ of subgraph folds (Notation 4.4) is called an *iterated subgraph fold*.

LEMMA 5.2. *Suppose that $(Z, X)$ is a labeled graph pair and that $\gamma \in \pi_1(Z)$. Then $\gamma \in N = N(Z, X)$ if and only if $\gamma \in \operatorname{Ker}((Q_Z)_\#)$ for some iterated subgraph fold $Q : (Z, X) \to (Z', X')$.*

*Proof.* Suppose that $q : (Z, X) \to (Z_1, X_1)$ is a subgraph fold. If $q_Z$ does not fold a bigon, then $q_Z$ is a homotopy equivalence. If $q_Z$ does fold a bigon then, up to homotopy equivalence, $q_Z$ is defined by adding a disk along a loop in $X$ with contractible image under $f_X$. It follows that, for an iterated subgraph fold $Q$, the map $Q_Z$ is defined, up to homotopy equivalence, by adding disks along loops in $X$ whose $f_X$-images are contractible. This implies $\operatorname{Ker}((Q_Z)_\#) \subset N$.

Conversely, suppose that $Q : (Z, X) \to (Z', X')$ is defined by iteratively folding subgraph edges until $f_{X'}$ is an immersion. Then $\operatorname{Ker}((f_X)_\#) \subset \operatorname{Ker}((Q_Z)_\#)$ and hence $N \subset \operatorname{Ker}((Q_Z)_\#)$. □

We will need the following simple observation in the proof of the next proposition.

*Remark* 5.3. If $Z$ is a wedge $A \vee B$ and if $q_Z : Z \to Z_1$ is a folding map that folds two edges in $A$, then the induced map from $A$ to $q_Z(A)$ is also a fold, the induced map from $B$ to $q_Z(B)$ is a homeomorphism, and $Z_1 = q_Z(A) \vee q_Z(B)$.

The next proposition is the relative version of Proposition 3.4. Before proving it, we state an immediate corollary that we will apply in proving the Main Proposition.

PROPOSITION 5.4. *Suppose that $(Z, X)$ is an invariant labeled graph pair for $H$. The tightening procedure produces from $(Z, X)$ a tight invariant labeled graph pair $(\hat{Z}, \hat{X})$ for $H$ with $\operatorname{rr}(\hat{Z}, \hat{X}) \leq \operatorname{rr}(Z, X)$. If $(f_Z)_\#$ is not injective, but $(f_X)_\#$ is injective, then $\operatorname{rr}(\hat{Z}, \hat{X}) < \operatorname{rr}(Z, X)$.*

COROLLARY 5.5. *If $(Z, X)$ is an invariant labeled graph pair for $H$ of minimal relative rank and if $(f_X)_\#$ is injective, then $(f_Z)_\#$ is injective.*

*Proof.* That $(\hat{Z}, \hat{X})$ is tight and invariant follows the definition of tightening. Lemma 4.7(4) implies that $\operatorname{rr}(\hat{Z}, \hat{X}) \leq \operatorname{rr}(Z, X)$.

It remains to prove the last sentence of the proposition. A single step in the tightening procedure is to fold $q : (Z_0, X_0) \to (Z_1, X_1)$ and then add a loop if necessary to form $(Z_2, X_2)$. Let $\operatorname{inc} : Z_1 \to Z_2$ the inclusion map, let $*_i$ be the basepoint of $Z_i$ and let $N_i = N(Z_i, X_i)$ for $i \in \{0, 1, 2\}$.

We will show that if $\operatorname{rr}(Z_2, X_2) = \operatorname{rr}(Z_0, X_0)$ and if $N_0 \neq \operatorname{Ker}((f_{Z_0})_\#)$ then $N_2 \neq \operatorname{Ker}((f_{Z_2})_\#)$. This will finish the proof of the proposition. Indeed, at the very beginning of the tightening process, $(Z_0, X_0) = (Z, X)$; in this case $N_0$ is trivial and $\operatorname{Ker}((f_{Z_0})_\#)$ is nontrivial. If $\operatorname{rr}(\hat{Z}, \hat{X}) = \operatorname{rr}(Z, X)$, then $\operatorname{rr}(Z_2, X_2) = \operatorname{rr}(Z_0, X_0)$ for each tightening step and we conclude, by induction on the number of steps, that $N(\hat{Z}, \hat{X}) \neq \operatorname{Ker}((f_{\hat{Z}})_\#)$. But, this contradicts the fact that $(\hat{Z}, \hat{X})$ is tight. We conclude that $\operatorname{rr}(\hat{Z}, \hat{X}) < \operatorname{rr}(Z, X)$ as desired.

Assume now that $\operatorname{rr}(\hat{Z}_2, \hat{X}_2) = \operatorname{rr}(Z_0, X_0)$ and that $N_0 \neq \operatorname{Ker}((f_{Z_0})_\#)$.

*Step* 1 (Subgraph folds). If $q$ is a subgraph fold then $(Z_1, X_1) = (Z_2, X_2)$ and Lemma 5.2 implies that if $\gamma$ is any nontrivial element of $\operatorname{Ker}((f_{Z_0})_\#) \setminus N_0$ then $(q_{Z_0})_\#(\gamma)$ is a nontrivial element of $\operatorname{Ker}((f_{Z_1})_\#) \setminus N_1$. We may therefore assume that $q$ is not a subgraph fold.

*Step* 2 (The case that $N_1$ is trivial). By definition of the tightening procedure, $f_{X_0}$ must be an immersion; thus $N_0$ is trivial. Since $\mathrm{rr}(Z_2, X_2) = \mathrm{rr}(Z_0, X_0)$, Lemma 4.7(5) implies that $q_{Z_0}$ does not fold a bigon. In particular, $(q_{Z_0})_\#$ is injective and $(\mathrm{inc})_\#(q_{Z_0})_\#\mathrm{Ker}((f_{Z_0})_\#)$ is a nontrivial subgroup of $\mathrm{Ker}((f_{Z_2})_\#)$. If $N_1$, and hence $N_2$, is trivial, then we are done. We may therefore assume that $N_1$ is nontrivial.

*Step* 3 (Making a jump). By Lemma 4.3, $q$ is an exceptional fold. Let $p_1$ and $p_2$ be the points in $X_0$ that are identified by $\check{q}_{X_0}$; let $p = \check{q}_{X_0}(p_1) = \check{q}_{X_0}(p_2)$ be their image in $X_1$; for $i = 1, 2$, let $d_i$ be the path in $X_0$ connecting $*_0$ to $p_i$ in the definition of $(Z_2, X_2)$; and let $d = \check{q}_{X_0}(d_1)\check{q}_{X_0}(d_2^{-1})$ be the resulting loop in $X_1$. Recall that $\delta \in \mathbb{F}$ is represented by $f_{X_1}(d)$. By Lemma 4.7(6), $\Phi(\delta) \notin (f_{X_1})_\#(\pi_1(Z_1))$. The link of $p$ in $X_1$ is the union $L_1 \cup L_2$ where $L_j$ is the link of $p_j$ in $X_0$ for $j = 1, 2$. Let $\rho$ be any closed based path in $X_1$. If $\rho$ crosses $p$, entering through $L_1$ and exiting through $L_2$, or vice-versa, then we say that $\rho$ *makes a jump* at that crossing of $p$. Write $\rho = \rho_1 \ldots \rho_m$ as a concatenation of subpaths where we have subdivided each time $\rho$ makes a jump at $p$. The paths $\check{q}_{X_0}(d_1)$ and $\check{q}_{X_0}(d_2)$ connect $*_1$ to $p$. By inserting $\check{q}_{X_0}(d_1^{-1})\cdot\check{q}_{X_0}(d_1)\cdot\check{q}_{X_0}(d_2^{-1})\cdot\check{q}_{X_0}(d_2) = \check{q}_{X_0}(d_1^{-1})\cdot d\cdot\check{q}_{X_0}(d_2)$ or its inverse between the $\rho_i$'s, we produce closed paths $\rho_i^*$ based at $*_1$ such that $\rho \simeq \rho_1^*\rho_2^*\ldots\rho_s^*$ and such that $\rho_i^*$ either makes no jumps at $p$ or equals $d$ or $d^{-1}$.

*Step* 4 (Factoring $\beta \in \mathrm{Ker}((f_{X_1})_\#)$). Since $N_1$ is nontrivial, we may choose a nontrivial element $\beta \in \mathrm{Ker}((f_{X_1})_\#)$. Applying the above decomposition to a closed based path representing $\beta$, we conclude that $\beta$ is the alternating product of nontrivial elements $\mu_i, \nu_j \in \pi_1(X_1)$ where $\mu_i = (\check{q}_{X_0})_\#(\mu_i')$ for some $\mu_i' \in \pi_1(X_0)$ and where $\nu_j$ is represented by a multiple of $d$ or $d^{-1}$. Since $N_0$ is trivial, at least one $\nu_j$ must appear in this product.

*Step* 5 (Applying invariance to produce $\gamma_2 \in \mathrm{Ker}((f_{Z_2})_\#)$). We next construct an element $\gamma_2 \in \pi_1(Z_2)$ such that $(f_{Z_2})_\#(\gamma_2) = \Phi((f_{X_1})_\#(\beta))$.

By definition, $Z_2$ is the wedge of $Z_1$ with $Y = X(\{\Phi(\delta)\})$; we will think of $Z_1$ and $Y$ as subgraphs of $Z_2$. For each $\nu_j$, there exists a nontrivial element $\tau_j \in \pi_1(Y) \subset \pi_1(Z_2)$ such that $(f_{Z_2})_\#(\tau_j) = \Phi((f_{X_1})_\#(\nu_j))$. By the invariance property for $(Z_0, X_0)$, for each $\mu_i$, there exists a nontrivial element $\sigma_i' \in \pi_1(Z_0)$ such that $(f_{Z_0})_\#(\sigma_i') = \Phi((f_{X_0})_\#(\mu_i')) = \Phi((f_{X_1})_\#(\mu_i))$. Let $\sigma_i = (q_{Z_0})_\#(\sigma_i') \in \pi_1(Z_1) \subset \pi_1(Z_2)$; then $(f_{Z_2})_\#(\sigma_i) = \Phi((f_{X_1})_\#(\mu_i))$. Define $\gamma_2$ to be the alternating product of the $\sigma_i$'s and $\tau_j$'s in the same order as $\beta$ is the product of the $\mu_i$'s and $\nu_j$'s. Then $(f_{Z_2})_\#(\gamma_2) = \Phi((f_{X_1})_\#(\beta))$ as promised. Since $\beta \in \mathrm{Ker}((f_{X_1})_\#)$, $\gamma_2 \in \mathrm{Ker}((f_{Z_2})_\#)$.

*Step* 6 ($\gamma_2 \notin N_2$). It suffices to show that $\gamma_2 \notin N_2$. Since $N_0$ is trivial and $\mu_i' \in \pi_1(X_0)$, $(f_{X_0})_\#(\mu_i')$ and hence $(f_{Z_2})_\#(\sigma_i) = \Phi((f_{X_0})_\#(\mu_i'))$ is nontrivial. Suppose that $Q : (Z_2, X_2) \to (Z', X')$ is an iterated subgraph fold. Denote $Q_{Z_2}(Z_1)$ by $Z_1'$. It follows from Remark 5.3 that $Q_{Z_2}$ induces a homeomorphism from $Y$ to $Y' = Q_{Z_2}(Y)$, that $Z'$ is the wedge of $Z_1'$ and $Y'$ at $Q_{Z_2}(*_2)$, and that $(Q_{Z_2})_\#(\tau_j)$ represents a nontrivial element of $\pi_1(Y')$. Since $\text{Ker}((Q_{Z_2})_\#) \subset \text{Ker}((f_{Z_2})_\#)$, each $(Q_{Z_2})_\#(\sigma_i)$ represents a nontrivial element of $\pi_1(Z_1')$. Since $(Q_{Z_2})_\#(\gamma_2)$ is an alternating concatenation of nontrivial elements in $\pi_1(Y')$ and in $\pi_1(Z_1')$, it represents a nontrivial element of $\pi_1(Z')$. Lemma 5.2 implies that $\gamma_2 \notin N_2$. $\square$

## 6. Conclusion

*Proof of Main Proposition.* Remark 4.2 and the fact that relative rank is a nonnegative integer imply that there is an invariant labeled graph pair for $H$ of minimal relative rank. Proposition 5.4 therefore implies that there is a tight invariant labeled graph pair $(Z, X)$ for $H$ of minimal relative rank. Let $T$ be a maximal tree for $Z$ containing $*$ such that $T \cap X$ is a maximal tree for $X$. Let $\{a_1, \cdots, a_m, b_1, \cdots, b_r\}$ and $\{a_1, \cdots, a_m\}$ be the resulting bases for $Z^\#$ and $X^\#$ respectively. Set $A = \{a_1, \cdots, a_m\}$ and $B = \{b_1, \cdots, b_r\}$. For $1 \leq j \leq m$, let $w_j$ be the reduced word in $A \cup B$ such that $\Phi(a_j) = w_j$, let $r_j = ta_jt^{-1}w_j^{-1}$, and let $C = \{r_1, \cdots, r_m\}$. It suffices to prove that $\langle t, A, B | C \rangle$ is a presentation for $H$.

Let $F$ denote the free group with basis $\{t\} \cup A \cup B$, and let $\eta : F \to H$ be the natural map.

Set $X_0 = X(A)$ and for each $i \geq 1$, set $X_i = X(A \cup \bigcup_{l=0}^{d-1} \Phi^l(B))$. Then each $(X_{i+1}, X_i)$ is an invariant labeled graph pair for $H$ of minimal relative rank and $(f_{X_0})_\#$ is injective. So, by induction using Corollary 5.5, $(f_{X_i})_\#$ is injective for all $i \geq 0$. Thus, if $S_d' = A \cup \bigcup_{l=0}^{d} t^l B t^{-l} \subset F$ then the restriction of $\eta$ to $\langle S_d' \rangle$ is injective.

It is clear that $C$ is contained in the kernel of $\eta$. It remains to show that if $k$ is in the kernel of $\eta$ then $k$ is in the normal closure $N$ of $C$ in $F$. After replacing $k$ by $t^z k t^{-z}$ for some $z \in \mathbb{Z}$, we may assume that $k$ factors into a product where each term has the form $t^d x t^{-d}$ for some $d \geq 0$ and for some $x \in A \cup B$. We now prove by induction on $d$ the claim that each such term can be written as $nv$ where $n \in N$ and $v \in \langle S_d' \rangle$.

This is trivial for $d = 0$ and clear for $d = 1$ since $tb_i t^{-1} \in S_1'$, $ta_j t^{-1} = r_j w_j$ and $w_j \in \langle S_0' \rangle$. Assume now that $d > 1$ and that the result holds for $d - 1$. Then $t^d x t^{-d} = t(t^{d-1} x t^{-(d-1)})t^{-1} = tn_1 v_1 t^{-1} = (tn_1 t^{-1})(tv_1 t^{-1}) = n_2(tv_1 t^{-1})$ where $n_1, n_2 \in N$ and $v_1 \in \langle S_{d-1}' \rangle$. Now $tv_1 t^{-1}$ can be written as a product where each term is either in $\langle S_d' \rangle$ or is some $ta_j t^{-1}$. Since terms of the latter type factor into an element of $N$ and an element of $\langle S'_0 \rangle$, $t^d x_j t^{-d}$ factors into a product where each term is either in $N$ or in $\langle S'_d \rangle$. Since $N$ is normal, we may assume that all the terms in $N$ precede all the terms in $\langle S'_d \rangle$. This verifies our claim.

We have shown that $k$ is a product where each term is either in $N$ or in some $\langle S'_d \rangle$. Pushing all the $N$-terms to the front, we have $k = nv$ where $n \in N$ and $v \in \langle S'_d \rangle$ for some $d \geq 0$. Since $\eta(k)$ and $\eta(n)$ are trivial in $H$, so is $\eta(v)$. The injectivity of the restriction of $\eta$ to $\langle S'_d \rangle$ implies that $v$ must be the trivial word in $S'_d$ and so $k \in N$ as desired. □

Rutgers University, Newark, NJ
*E-mail address*: feighn@andromeda.rutgers.edu

Herbert H. Lehman College, CUNY, Bronx, NY
*E-mail address*: michael@alpha.lehman.cuny.edu